\documentclass[preprint]{elsarticle}       

\usepackage{lineno,hyperref}
\usepackage{rotating}
\usepackage{makeidx}
\usepackage{float}
\usepackage{epsfig}
\usepackage{amsmath,amssymb}
\usepackage[latin1]{inputenc}
\usepackage{epic,eepic}
\usepackage{overpic}
\usepackage{color}
\usepackage{amsfonts}
\usepackage{graphicx}
\usepackage{mwe} 
\usepackage{subcaption}

\usepackage{algorithm}
\usepackage{algorithmicx}
\usepackage{algpseudocode}
\usepackage{enumitem}

\newenvironment{demo}{\smallskip\noindent{{\it Proof.}}\hskip \labelsep}%
            {\hfill\penalty10000\raisebox{-.09em}{\large\bf\rm $\blacksquare$}\par\medskip}

\newtheorem{theorem}{Theorem}[section]

\newtheorem{remark}{Remark}[section]

\DeclareMathOperator*{\argmin}{arg\,min}

\def\xx{\mathbf{x}}
\journal{}

\begin{document}

\small

\begin{frontmatter}

\title{Weighted Essentially Non-Oscillatory Shepard method}\tnotetext[label1]{The  fourth author has been supported through project CIAICO/2021/227 (Proyecto financiado por la Conselleria de Innovaci\'on, Universidades, Ciencia y Sociedad digital de la Generalitat Valenciana), by grant PID2020-117211GB-I00 and by PID2023-146836NB-I00 funded by MCIN/AEI/10.13039/501100011033.}

\author[TAU]{David Levin}
\ead{levindd@gmail.com}
\author[UV]{Jos\'e M. Ram\'on}
\ead{Jose.Manuel.Ramon@uv.es}
\author[UPCT]{Juan Ruiz-\'Alvarez}
\ead{juan.ruiz@upct.es}
\author[UV]{Dionisio F. Y\'a\~nez}
\ead{Dionisio.Yanez@uv.es}

\date{Received: date / Accepted: date}

\address[TAU]{School of Mathematical Sciences. Tel-Aviv University, Tel-Aviv (Israel).}
\address[UV]{Departamento de Matem\'aticas. Universidad de Valencia, Valencia (Spain).}
\address[UPCT]{Departamento de Matem\'atica Aplicada y Estad\'istica. Universidad  Polit\'ecnica de Cartagena, Cartagena (Spain).}


\begin{abstract}

Shepard method is a fast algorithm that has been classically used to interpolate
scattered data in several dimensions. This is an important and well-known
technique in numerical analysis founded in the main idea that data that is far
away from the approximation point should contribute less to the resulting
approximation. Approximating piecewise smooth functions in $\mathbb{R}^n$ near
discontinuities along a hypersurface in $\mathbb{R}^{n-1}$ is challenging for the
Shepard method or any other linear technique for sparse data due to the inherent
difficulty in accurately capturing sharp transitions and avoiding oscillations.
This letter is devoted to constructing a non-linear Shepard method using the
basic ideas that arise from the weighted essentially non-oscillatory
interpolation method (WENO). The proposed method aims to enhance the accuracy and
stability of the traditional Shepard method by incorporating WENO's adaptive and
nonlinear weighting mechanism. To address this challenge, we will nonlinearly
modify the weight function in a general Shepard method, considering any weight
function, rather than relying solely on the inverse of the distance
squared. This approach effectively reduces oscillations near discontinuities and improves the overall interpolation quality. Numerical experiments demonstrate the superior performance of the new method in handling complex datasets, making it a valuable tool for various applications in scientific computing and data analysis.

\end{abstract}
\end{frontmatter}
\vspace{-0.5cm}

\section{Introduction}

The Shepard method, introduced by Donald Shepard in 1968 \cite{shepard1968}, is a widely used algorithm for interpolating scattered data in multiple dimensions. Traditionally, the weight function in the Shepard method is the inverse of the squared distance. However, this weight function can be generalized to any function, such as Gaussian \cite{gaussian}, inverse multiquadratic \cite{inverse_multiquadratic}, Mat\'ern \cite{matern}, and Wendland functions \cite{wendland,wendland2002}. These functions must have compact support and form a partition of unity to ensure the interpolation's accuracy and stability \cite{fasshauer2007}.


The linear Shepard method performs well with smooth functions, such as the Franke's function. However, it tends to produce non desirable diffusion effects near discontinuities, which can significantly degrade the approximation quality. In this letter, we propose a nonlinear Shepard's method that aims to mitigate this diffusion. By modifying the weight function nonlinearly, we achieve results near discontinuities that are comparable to the Shepard's method performance in smooth regions.

The letter is organized as follows: The second section discusses the construction of the nonlinear method. We prove some properties in Section 3. Finally, Section \ref{NE} presents a series of numerical experiments comparing the performance of the linear and nonlinear methods and checking the theoretical results.

\section{Construction of the non-linear method}

Given an open and bounded domain $\Omega\subseteq \mathbb{R}^n$, a set $\chi_N=\{\xx_i\in \Omega:i=1,\hdots,N\}$ of $N$ distinct nodes and a corresponding set of function values $\mathcal{F}_N=\{f_i=f(\xx_i):i=1,\hdots,N\}$, where $f:\Omega\to\mathbb{R}$. We assume the points $\chi_N$ are quasi-uniformly distributed in $\Omega$, with a fill distance $h =\sup_{x\in\Omega}\min_{\xx_i\in\chi_N}\|\xx-\xx_i\|$. Shepard's method can be viewed as a special case of the moving least squares technique (see \cite{davidlevin}), where the degree of the polynomial is zero, i.e., a constant. Consequently, if we consider a non-negative, compactly supported, radial function $\omega_i:\Omega\to\mathbb{R}$, and define
$\omega_i(\xx)=\omega(\frac{\|\xx_i-\xx\|}{h}),$ the Shepard's interpolant is defined as:
$\mathcal{I}_{\text{S}}(\xx)=\sum_{i=1}^N \frac{\omega_i(\xx)}{\sum_{j=1}^N\omega_j(\xx)}f_i=\sum_{i=1}^N W_i(\xx)f_i,$
where $W_i:\Omega\to \mathbb{R}$ are defined as $W_i(\xx)=\frac{\omega_i(\xx)}{\sum_{j=1}^N\omega_j(\xx)}$. It follows that for $i=1,\hdots,N$, we have
$ 0\leq W_i(\xx)\leq 1,\, \text{and} \,\sum_{j=1}^NW_j(\xx)=1, \forall \xx\in\mathbb{R}^n$.
We call the $W_i$ ``{\it optimal Shepard weights}''.
\begin{remark}
In Shepard's method, the requirement that the function
$\omega(r)$ is compactly supported can be relaxed to a condition of rapid decay as $r\to\infty$.
Interpolation is achieved if
$\omega(0)=\infty$, (see \cite{davidlevin}).
\end{remark}

As described above, Shepard's method is based on the distances between the evaluation point and the given nodes. However, if a discontinuity contaminates the data, or they present a strong gradient, It is then recommended to use only nodes that are free of singularities. To achieve this, the optimal Shepard weights are replaced with non-linear ones, defined as follows:
$\mathcal{W}_i(\xx)=\frac{\alpha_i(\xx)}{\sum_{j=1}^N\alpha_j(\xx)},\quad \alpha_i(\xx)=\frac{W_i(\xx)}{(\epsilon+I_i)^t}.$
 The values $\{I_i\}$ are smoothness indicators, see \cite{doi:10.1137/070679065}, which identify whether a node is near a discontinuity and thus should have a negligible contribution to the approximation. $\epsilon$ and $t$ are parameters chosen to achieve maximal approximation order. In our case, we take $\epsilon=10^{-14}$ and $t=4$. In the next subsection, we present a way to construct the smoothness indicators.  With these ingredients the new WENO-Shepard's method is:
\begin{equation}\label{nolinear}
\mathcal{I}_{\text{WENO-S}}(\xx)=\sum_{i=1}^N \mathcal{W}_i(\xx)f_i.
\end{equation}
\subsection{The smoothness indicators}

We start dividing our domain $\Omega$ in some subdomains $\mathcal{S}_i$ with $\Omega\subset \cup_{i=1}^N\mathcal{S}_i$, which satisfy that
$\xx_i\in\mathcal{S}_i, \, i=1,\hdots,N$. In our case, we define the {\it stencils} as the balls centered in the data points, $\xx_i$ with a determined radius, $\delta_i$ i.e.,
\begin{equation}\label{stencilSi}
\mathcal{S}_i=\chi_N\cap B(\xx_i,\delta_i)=\{\xx_j\in\Omega:\|\xx_j-\xx_i\|<\delta_i\}=\{\xx_{j_i}:j_i=1,\hdots,N_i\}.
 \end{equation}

 To determine the radius, we rely on the general results presented in \cite{fasshauer2007, davidlevin, wendland2002},  which assume that $\delta_i=ch$ where $c$ is a constant and $h$ is the fill distance.

To design the smoothness indicators, we aim for two properties indicated, for example, in \cite{ABM}:
\begin{enumerate}[label={\bfseries P\arabic*}]
\item\label{P1sm1d} The order of a smoothness indicator that is free of discontinuities is $h^2$, i.e.
$I_{i}=O(h^2) \,\, \text{if}\,\, f \,\, \text{is smooth in } \,\, \mathcal{S}_{i}.$
\item\label{P2sm1d} When a discontinuity crosses the stencil $\mathcal{S}_{i}$ then
$I_{i} \nrightarrow 0 \,\, \text{as}\,\, h\to 0.$
\end{enumerate}

Thus, we solve the linear least square problem for each $i=1,\hdots,N,$
$p_i=\argmin_{p\in\Pi_1(\mathbb{R}^2)} \sum_{j_i=1}^{N_i} (f(\xx_{j_i})-p(\xx_{j_i}))^2,$
and define
\begin{equation}\label{indicador1}
I_i=\frac{1}{N_i}\sum_{j_i=1}^{N_i}|f(\xx_{j_i})-p_i(\xx_{j_i})|.
\end{equation}
It is clear that the $I_i$ defined in \eqref{indicador1} satisfy \ref{P1sm1d} and \ref{P2sm1d}. With these definitions we can prove some properties about the smoothness, the order of accuracy, and the behaviour close to the discontinuities presented by the approximant.

\section{Properties of the new method}\label{proper}

In this section, we prove some characteristics about our new method. Principally, we want to conserve the properties of the linear one, Shepard's method,
and to improve the approximation in the zones close to the discontinuities. Let us start with the smoothness of the non-linear method. In this case, from the expression of the smoothness indicator,
Eq. \eqref{indicador1}, we can see that it is a number independent of the variable $x$. Then, we only need the smoothness of the weight function $\omega$. This is similar to the linear Shepard's operator. 
We can summarize this property with the following result (Th. \ref{teo1}).

\begin{theorem}\label{teo1}
Let $\nu\in\mathbb{N}$, $\Omega\subset \mathbb{R}^n$, $\omega$ be a function with $\omega\in\mathcal{C}^\nu(\Omega)$,
and let $\mathcal{I}_{\text{\it WENO-S}}$ be the WENO Shepard's approximation defined
in Eq. \eqref{nolinear}, then $\mathcal{I}_{\text{\it WENO-S}}\in \mathcal{C}^\nu(\Omega)$.
\end{theorem}

Regarding the order of accuracy in the smooth regions, it is evident that the new operator reproduces constants, as it satisfies the condition that:
$\sum_{i=1}^N \mathcal{W}_i(\xx)=1,\,\, \forall\, \xx \in \mathbb{R}^n.$
As $\omega$ is compactly supported, then we have the next theorem.

\begin{theorem}\label{teo2}
Let $\Omega\subset \mathbb{R}^n$. If $f\in\mathcal{C}^1(\overline{\Omega})$, $\chi_N=\{\xx_i\in \Omega:i=1,\hdots,N\}$ are quasi-uniformly distributed with fill distance $h$, the weight function $\omega$ is compactly supported with support size $c$, then
$$\|f-\mathcal{I}_{\text{WENO-S}}\|_\infty\leq Ch\max_{\boldsymbol{\xi}\in \overline{\Omega}}|D^{\boldsymbol{\alpha}} f(\boldsymbol{\xi})|, \quad |\boldsymbol{\alpha}|=1,$$
where $C$ is a constant independent of $h$.
\end{theorem}

Finally, we analyze the behaviour of the operator $\mathcal{I}_{\text{WENO-S}}$ close to the discontinuities.
For this purpose, we consider a $(n-1)$-hypersurface $\Gamma$ defined by a function $\gamma:\mathbb{R}^n\to \mathbb{R}$ being
$\Omega^+=\{\xx\in \Omega: \gamma(\xx)\geq 0\}, \quad \Omega^-=\Omega\setminus \Omega^+,$
and we suppose that
\begin{equation}\label{tildef}
\tilde{f}(\xx)=\left\{
              \begin{array}{ll}
                f_1(\xx), & \xx\in \Omega^+,\\
                f_2(\xx), & \xx\in \Omega^-,
              \end{array}
            \right.
\end{equation}
with $f_1\in \mathcal{C}^1(\overline{\Omega^+})$ and $f_2\in \mathcal{C}^1(\overline{\Omega^-})$.

\begin{remark}{\bf Diffusion in Shepard's method}
When applying the linear Shepard method to discontinuous data with a radial weight function $\omega$ of support $c$, the resulting approximation exhibits a significant diffusion at all points within a distance $\le ch$ from the discontinuity $\Gamma$. In the numerical experiments in Section \ref{NE} we use Wendland's compactly supported weight functions (see \cite{wendland95}).
\end{remark}

The following Theorem shows that the diffusion area is significantly reduced when using the non-linear WENO-Shepard's method.
We use the above definitions of $\Gamma$, $\Omega^+$ and $\Omega^-$ and of $\tilde{f}$ defined in Eq. \eqref{tildef}.
We consider a weight function $\omega$ of support size $c$, $c\ge 2$, and such that $\omega(c-\epsilon_0)\ge C_0>0$. We further assume that $\Gamma$ is a smooth hypersurface.

\begin{theorem}\label{teo3new}
Let $\xx_0\in\Omega^+$ be at distance $h(1+\epsilon_0)$ from $\Gamma$. Then
$|\tilde{f}(\xx_0)- \mathcal{I}_{\text{WENO-S}}(\xx_0)|=O(h).$
\end{theorem}

\begin{demo}
Define the support of data points respect to $\xx_0$ as $S(\xx_0)=\{\xx_i\in\chi_N:\omega_i(\xx_0)>0\}$. It turns out that for any point $\xx_i\in S(\xx_0)\cap \Omega^-$ the smoothness indicator $I_i$ is $O(1)$. Furthermore, there exist points $\xx_i\in S(\xx_0)\cap \Omega^+$ such that the smoothness indicator $I_i$ is $O(h^2)$ as $h\to 0$. To prove this we examine the ball of radius $ch$ centered at $\xx_0$, $B(\xx_0, ch)$. We also consider all the points in $\Omega^+$ that are at least at a distance $ch$ away from $\Gamma$, and denote this set as $D^+$. Since $\xx_0\in \Omega^+$, and is located at a distance $h(1+\epsilon_0)$ from $\Gamma$, and given that $\Gamma$ is smooth, it follows that for sufficiently small $h$, the intersection $D^+\cap B(\xx_0, ch)$ contains a ball of diameter $h(1+\epsilon_0)$. Such a ball must contain a point $\xx_i\in S(\xx_0)\cap \Omega^+$, and for such a point
 $I_i=O(h^2)$ as $h\to 0$.

We define $\mathcal{K}=\{i : I_i=O(h^2), 1\leq i\leq N\}$.

It follows that
$$
\alpha_i(\xx)=\left\{
           \begin{array}{ll}
             O(1), & i\in \mathcal{K},\\
             O(h^{-2t}), & i\notin \mathcal{K},
           \end{array}
         \right.
$$
therefore $\sum_{i=1}^N \alpha_i(\xx_0)=O(h^{-2t})$ and $\mathcal{W}_i(\xx_0)=O(h^{2t})$ for $i\notin \mathcal{K}$. Then, we get
\begin{equation*}
\begin{split}
\tilde{f}(\xx_0)- \mathcal{I}_{\text{WENO-S}}(\xx_0)&=\sum_{i=1}^N \mathcal{W}_i(\xx_0) (\tilde{f}(\xx_0)-\tilde{f}_i)=\sum_{i\in\mathcal{K}} \mathcal{W}_i(\xx_0) (f_1(\xx_0)-f_1(\xx_i))+\sum_{i\notin\mathcal{K}} \mathcal{W}_i(\xx_0) (\tilde{f}(\xx_0)-\tilde{f}_i)\\
&=\sum_{i\in\mathcal{K}} \mathcal{W}_i(\xx_0) O(h)+\sum_{i\notin\mathcal{K}} O(h^{2t}) (\tilde{f}(\xx_0)-\tilde{f}_i)=O(h).
\end{split}
\end{equation*}

The same result holds if $\xx_0\in\Omega^-$ is at distance $h(1+\epsilon_0)$ from $\Gamma$.

\end{demo}
\section{Numerical experiments}\label{NE}
In this section, we perform some numerical tests to check the theoretical results proved in last section. We divide the section in two parts: firstly we analyze
the order of accuracy in the smooth zones.  In the second part, we will use a function with a jump discontinuity determined by a curve. All the examples are designed in $[0,1]^2\subset \Omega \subseteq \mathbb{R}^2$.

\subsection{Order of accuracy}
In order to study the behaviour in the smooth zones, we will approximate the Franke's function
using as nodes the regular grid $\chi_{2^l+1}=\{(i/2^l,j/2^l):i,j=0,\hdots,2^l\}$, and a set of $(2^l+1)^2$ Halton scattered data \cite{halton}. We denote the fill distance as $h_l$ and the errors as $e^l_i=|f(\mathbf{z}_i)-\mathcal{I}^l(\mathbf{z}_i)|$ with $\{\mathbf{z}_i:1\leq i\leq N_{\text{eval}}\}$ the set of points where we approximate the function. Finally we denote the maximum, discrete $\ell^2$ norms and rates as:
\begin{equation*}
\begin{split}
\text{MAE}_l=\max_{i=1,\hdots,N_{\text{eval}}}e^l_i, \quad \text{RMSE}_l=\left(\frac{1}{N_{\text{eval}}}\sum_{i=1}^{N_{\text{eval}}}(e_i^l)^2\right)^{\frac12}, \quad r_l^{\infty}=\frac{\log(\text{MAE}_{l-1}/\text{MAE}_{l})}{\log(h_{l-1}/h_{l})}, \quad r_l^2=\frac{\log(\text{RMSE}_{l-1}/\text{RMSE}_{l})}{\log(h_{l-1}/h_{l})}.
\end{split}
\end{equation*}
In these examples, we use the $\mathcal{C}^2$ and $\mathcal{C}^4$ Wendland's compactly supported functions (see \cite{wendland95} or chapter 11 in  \cite{fasshauer2007}) defined as $\omega_{\text{W2}}(r)=(1-\varepsilon r)^4_+(4\varepsilon r+1)$ and $\omega_{\text{W4}}(r)=(1-\varepsilon r)^6_+(35(r\varepsilon)^2 +18\varepsilon r+3)$ where $r$ is the Euclidean distance and $\varepsilon$ is the shape parameter. In our experiments, we take $\varepsilon=\lfloor \frac{2^l+1}{2}\rfloor/\sqrt{2}$. In Tables \ref{tabla3} and \ref{tabla4} we can see that the numerical order of accuracy is the one expected in the linear and non-linear methods (and sometimes improved by the last one). The results are very similar, therefore the behaviour of the new algorithm in the smooth zones is analogous to Shepard's method. In terms of the norm of the error, the linear method obtains a very slight advantage for this experiment.

\begin{table}[!ht]
\begin{center}
\begin{tabular}{lrrrrrrrrrrr}
(W2) & \multicolumn{2}{c}{Shepard} & &   \multicolumn{2}{c}{WENO-Shepard}& &\multicolumn{2}{c}{Shepard} & &   \multicolumn{2}{c}{WENO-Shepard}\\ \cline{1-3} \cline{5-6} \cline{8-9} \cline{11-12}
$l$ & $\text{MAE}_l$ & $r_l^{\infty}$ & &$\text{MAE}_l$ & $r_l^{\infty}$ & &$\text{RMSE}_l$ & $r_l^2$ & &$\text{RMSE}_l$ & $r_l^2$       \\
\hline
$4$   & 6.1891e-02 &        -&&  1.4160e-01 &       -    &&    1.5976e-02  &      - &&  4.4803e-02  &             \\
$5$   & 2.1657e-02 &  1.5149 &&  5.4509e-02 &  1.3773    &&    4.7667e-03  & 1.7448 &&  1.5903e-02  & 1.4943      \\
$6$   & 1.1315e-02 &  0.9365 &&  1.7315e-02 &  1.6545    &&    1.5991e-03  & 1.5758 &&  4.6052e-03  & 1.7879      \\
$7$   & 5.7795e-03 &  0.9693 &&  4.6431e-03 &  1.8989    &&    6.6941e-04  & 1.2563 &&  9.5190e-04  & 2.2744      \\
\hline \hline
(W4) & \multicolumn{2}{c}{Shepard} & &   \multicolumn{2}{c}{WENO-Shepard}& &\multicolumn{2}{c}{Shepard} & &   \multicolumn{2}{c}{WENO-Shepard}\\ \cline{1-3} \cline{5-6} \cline{8-9} \cline{11-12}
$l$ & $\text{MAE}_l$ & $r_l^{\infty}$ & &$\text{MAE}_l$ & $r_l^{\infty}$ & &$\text{RMSE}_l$ & $r_l^2$ & &$\text{RMSE}_l$ & $r_l^2$       \\
\hline
$4$   &    4.7978e-02 &         && 1.1916e-01  &          &&       1.2350e-02 &       - &&  3.6371e-02 &               \\
$5$   &    1.6849e-02 &  1.5097 && 4.3702e-02  & 1.4471  &&       3.6356e-03 &  1.7643 &&  1.2588e-02 &  1.5307      \\
$6$   &    8.8017e-03 &  0.9368 && 1.3690e-02  & 1.6746  &&       1.2311e-03 &  1.5623 &&  3.6177e-03 &  1.7989      \\
$7$   &    4.4875e-03 &  0.9719 && 3.6360e-03  & 1.9126  &&       5.1835e-04 &  1.2479 &&  7.3183e-04 &  2.3055      \\
\hline
\end{tabular}
\end{center}
\caption{Errors and rates using Shepard and WENO-Shepard methods for Franke's test function in a regular grid.}\label{tabla3}
\end{table}

\begin{table}[!ht]
\begin{center}
\begin{tabular}{lrrrrrrrrrrr}
(W2)& \multicolumn{2}{c}{Shepard} & &   \multicolumn{2}{c}{WENO-Shepard}& &\multicolumn{2}{c}{Shepard} & &   \multicolumn{2}{c}{WENO-Shepard}\\ \cline{1-3} \cline{5-6} \cline{8-9} \cline{11-12}
$l$ & $\text{MAE}_l$ & $r_l^{\infty}$ & &$\text{MAE}_l$ & $r_l^{\infty}$ & &$\text{RMSE}_l$ & $r_l^2$ & &$\text{RMSE}_l$ & $r_l^2$       \\
\hline
$4$  &        1.0652e-01  &       -  & &  1.5273e-01   &        -    & &        2.0217e-02   &    -    & &  4.8709e-02 &                 \\
$5$  &        6.2913e-02  &  0.6926  & &  6.5045e-02   & 1.1227  & &        8.1166e-03   & 1.2003  & &  1.7523e-02 &  1.3446         \\
$6$  &        2.9851e-02  &  1.6150  & &  2.6824e-02   & 1.9188  & &        3.5854e-03   & 1.7699  & &  6.0397e-03 &  2.3074         \\
$7$  &        1.3749e-02  &  1.0304  & &  1.4274e-02   & 0.8385  & &        1.8125e-03   & 0.9067  & &  2.1831e-03 &  1.3525         \\
\hline
\hline
(W4)& \multicolumn{2}{c}{Shepard} & &   \multicolumn{2}{c}{WENO-Shepard}& &\multicolumn{2}{c}{Shepard} & &   \multicolumn{2}{c}{WENO-Shepard}\\ \cline{1-3} \cline{5-6} \cline{8-9} \cline{11-12}
$l$ & $\text{MAE}_l$ & $r_l^{\infty}$ & &$\text{MAE}_l$ & $r_l^{\infty}$ & &$\text{RMSE}_l$ & $r_l^2$ & &$\text{RMSE}_l$ & $r_l^2$       \\
\hline
$4$   &    1.0459e-01 &       -  & & 1.3537e-01 &          &&     1.7996e-02 &     -      & & 4.0106e-02  &              \\
$5$   &    6.2727e-02 &  0.6724  & & 6.1287e-02 &  1.0422  &&     7.5610e-03 &  1.1405    & & 1.4581e-02  & 1.3308      \\
$6$   &    3.0973e-02 &  1.5286  & & 2.6223e-02 &  1.8390  &&     3.5886e-03 &  1.6144    & & 5.3386e-03  & 2.1765      \\
$7$   &    1.3367e-02 &  1.1169  & & 1.3816e-02 &  0.8517  &&     1.8699e-03 &  0.8664    & & 2.1466e-03  & 1.2109      \\
\hline
\end{tabular}
\end{center}
\caption{Errors and rates using Shepard and WENO-Shepard methods for Franke's test function evaluated at Halton points.}\label{tabla4}
\end{table}

\subsection{Approximation of piecewise smooth functions}

\begin{figure}[H]
\begin{center}
		\begin{tabular}{ccc}
\hspace{-1cm}	Shepard& \hspace{-1cm}	   WENO-Shepard & \hspace{-1cm}	   Original \\
\hspace{-1cm}		\includegraphics[width=6.3cm]{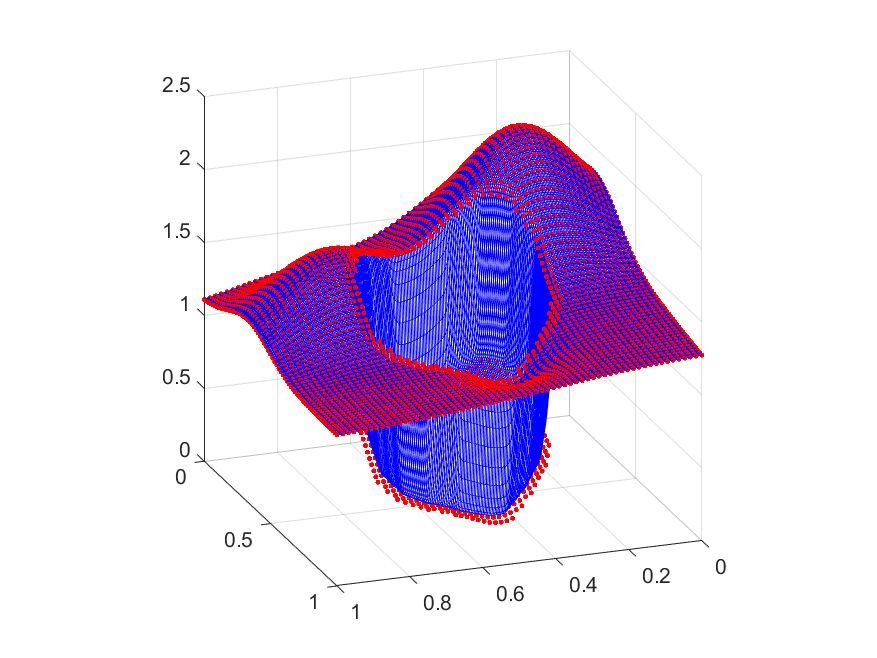} & \hspace{-1cm}	
			\includegraphics[width=6.3cm]{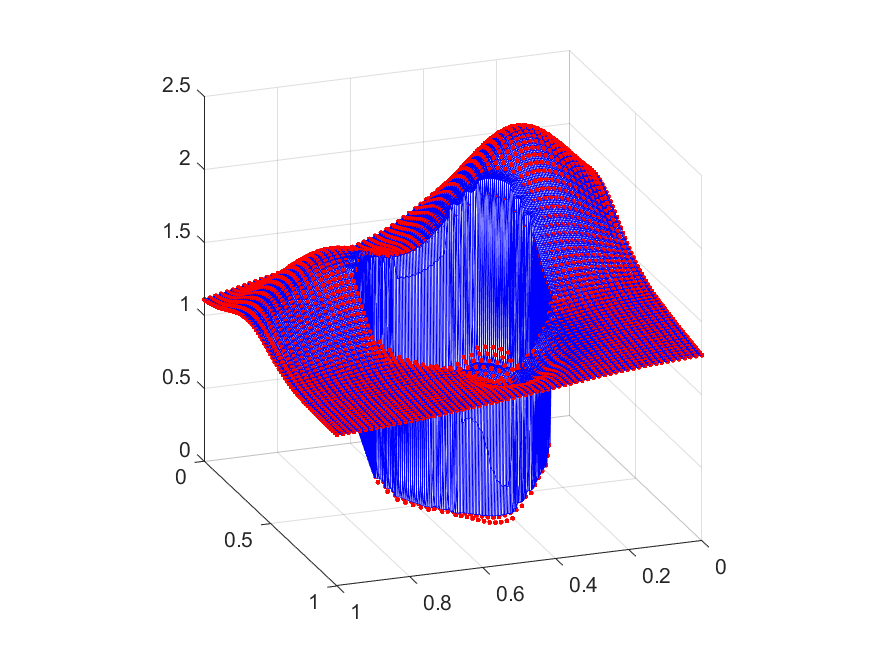} & \hspace{-1cm}	
			\includegraphics[width=6.3cm]{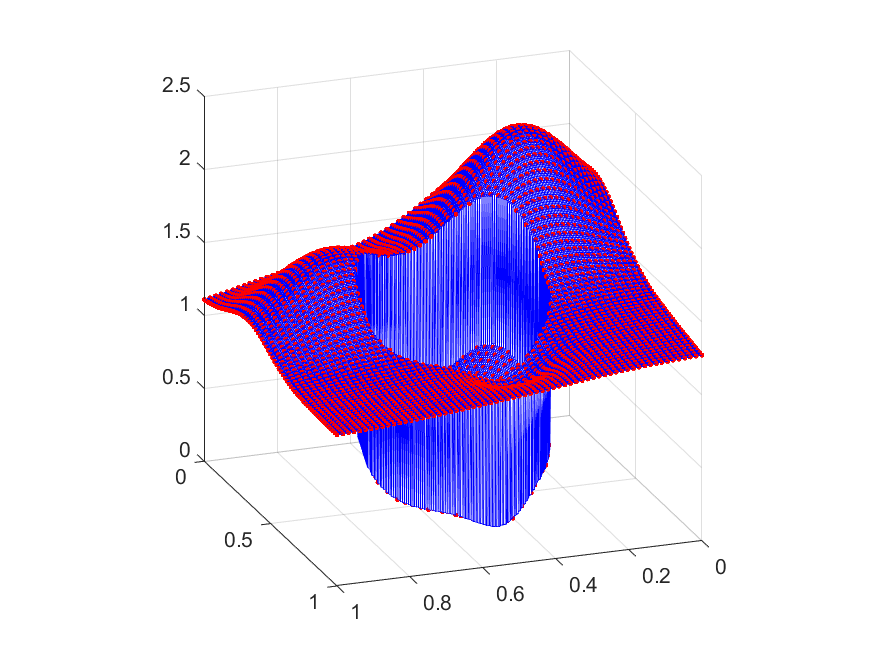}\\
\hspace{-1cm}		\includegraphics[width=6.3cm]{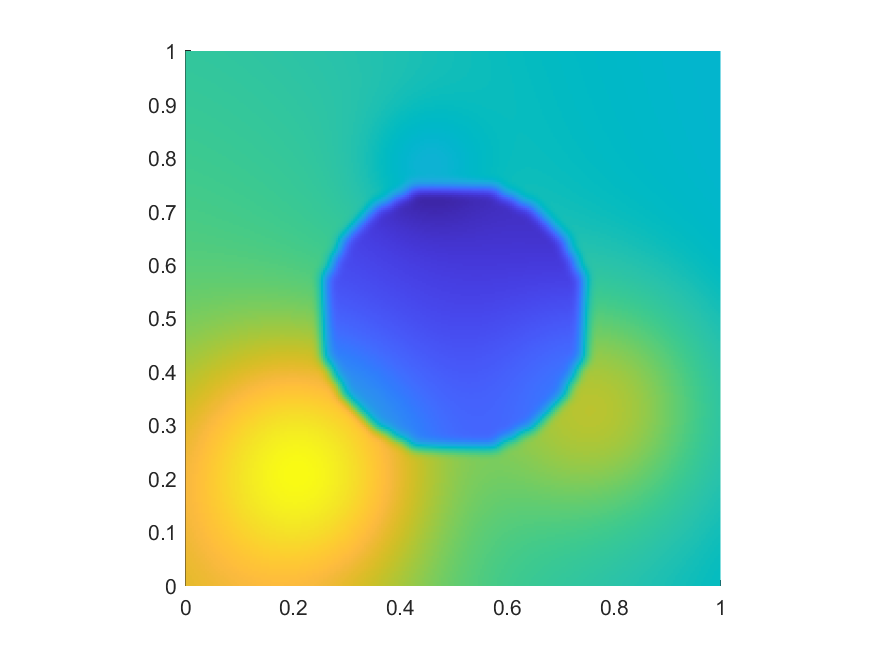} & \hspace{-1cm}	
			\includegraphics[width=6.3cm]{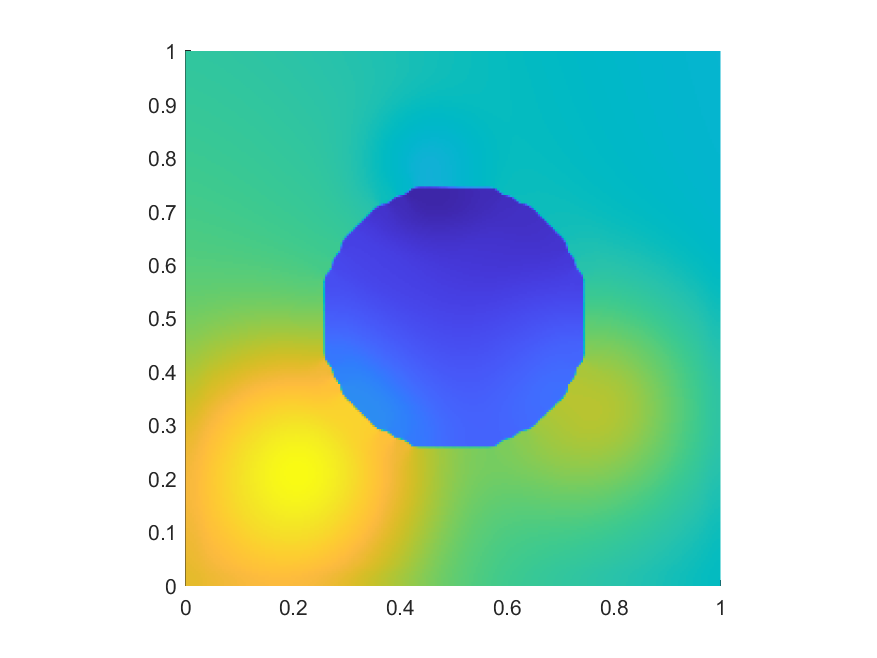} & \hspace{-1cm}	
			\includegraphics[width=6.3cm]{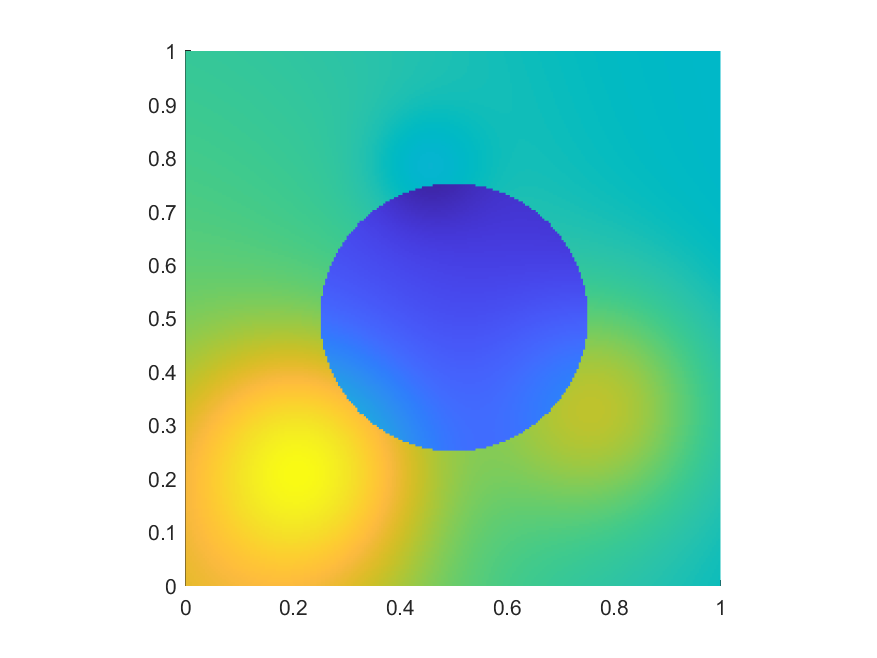}\\
\hspace{-1cm}		\includegraphics[width=6.3cm]{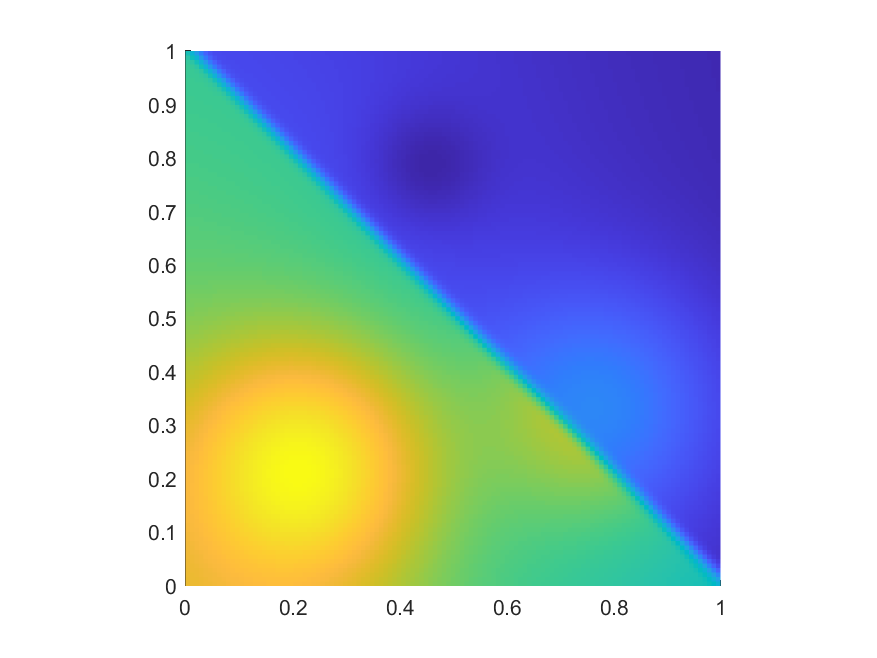} & \hspace{-1cm}	
			\includegraphics[width=6.3cm]{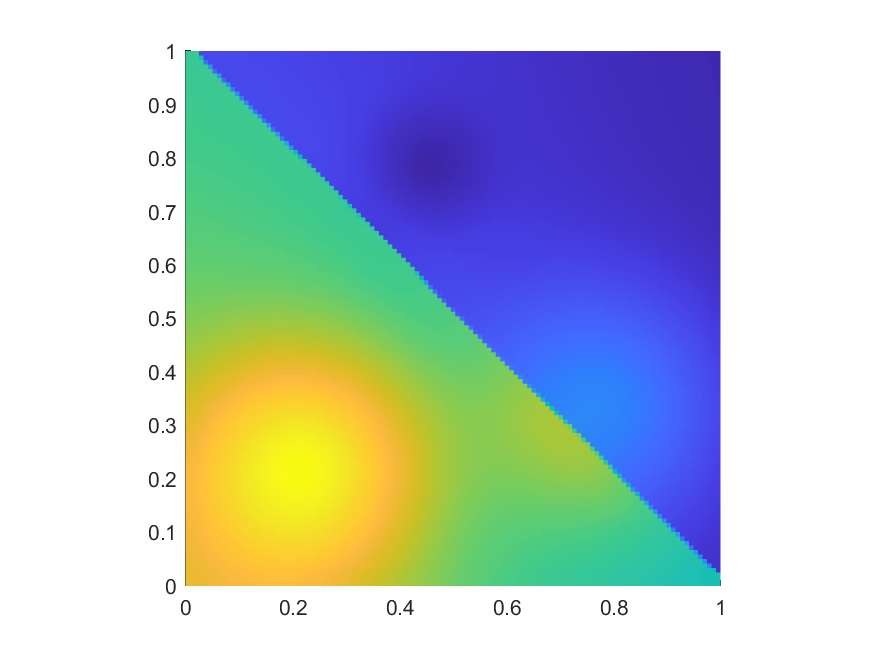} & \hspace{-1cm}
			\includegraphics[width=6.3cm]{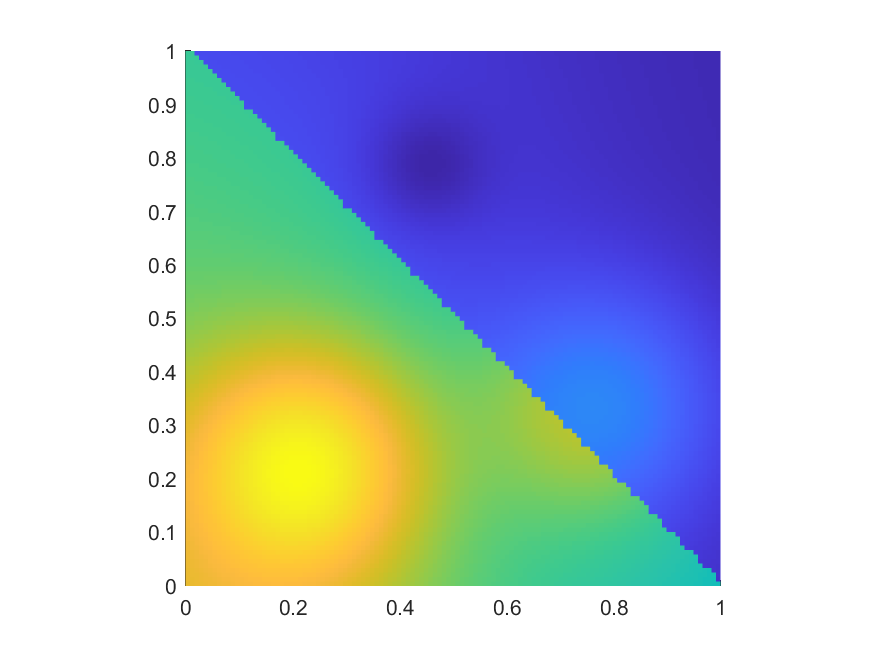}\\
		\end{tabular}
\end{center}
			\caption{Approximation to function $\tilde{f}$, Eq. \eqref{funciontildef}. In first line, a 3-D plot with a circumference as discontinuity curve and the next two rows in a 2-D plot with different discontinuity curves using regular grid data points.}
		\label{gridfiguracirc}
	\end{figure}

\begin{figure}[H]
\begin{center}
		\begin{tabular}{ccc}
\hspace{-1cm}	Shepard& \hspace{-1cm}	   WENO-Shepard & \hspace{-2cm}	Original\\
\hspace{-1cm}		\includegraphics[width=6.3cm]{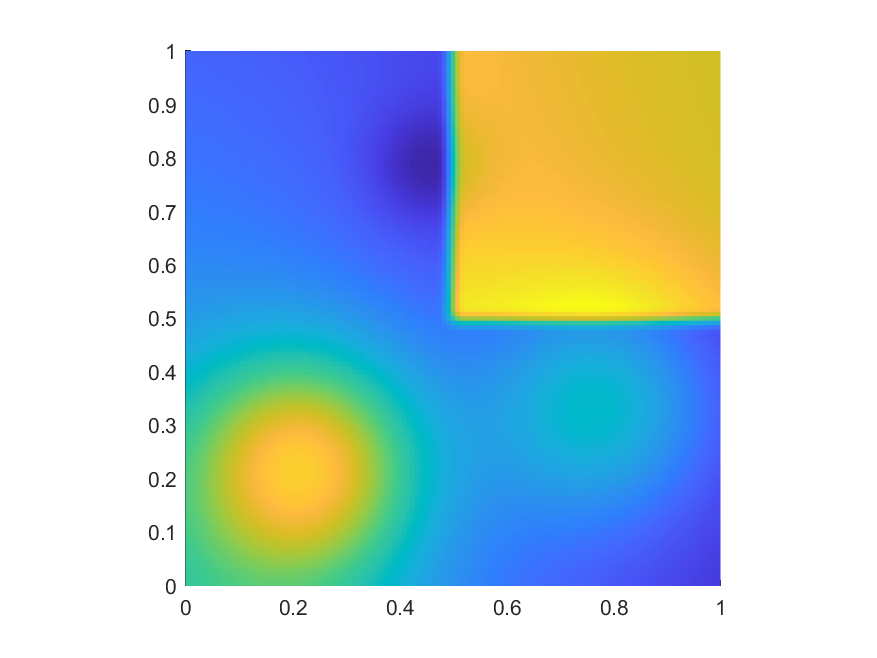} & \hspace{-1cm}	
			\includegraphics[width=6.3cm]{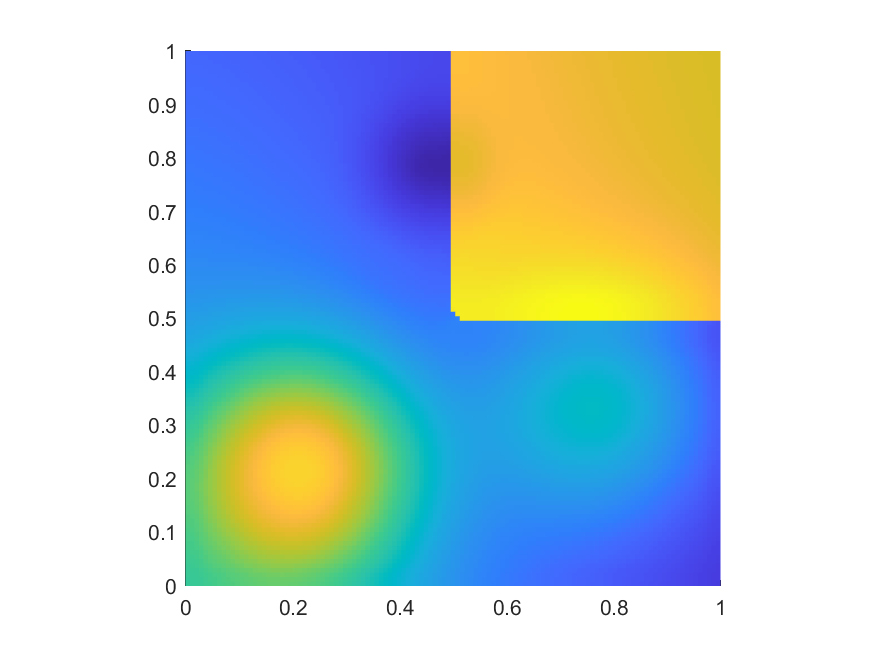} & \hspace{-1cm}
			\includegraphics[width=6.3cm]{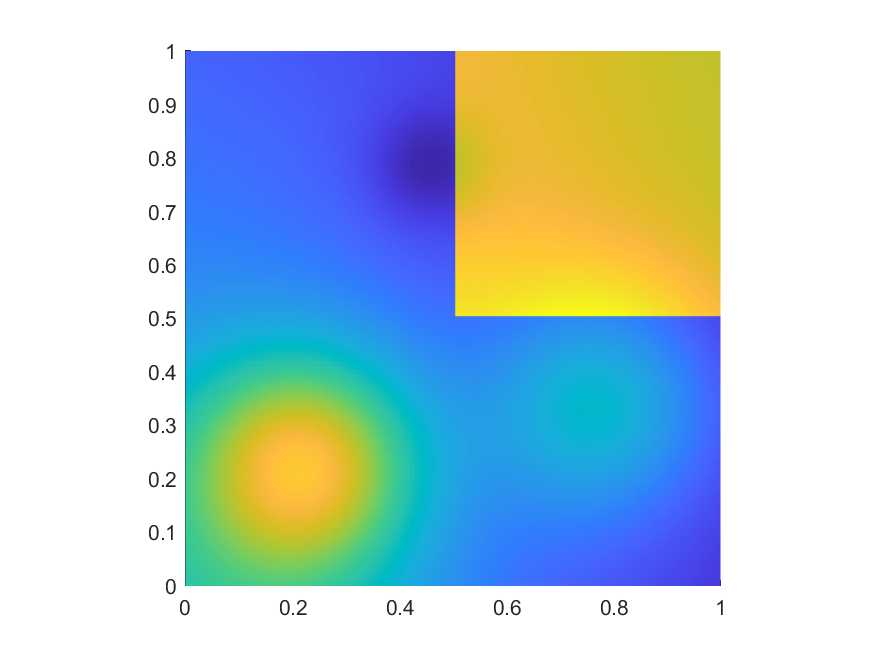}
		\end{tabular}
\end{center}
				\caption{Approximation to function $\tilde{f}$, Eq. \eqref{funciontildef}, in a 2-D plot using regular grid data points.}
		\label{gridfigura}
	\end{figure}

\begin{figure}[H]
\begin{center}
		\begin{tabular}{ccc}
\hspace{-1cm}	Shepard& \hspace{-1cm}	   WENO-Shepard & \hspace{-1cm}	Original\\
		\hspace{-1cm}		\includegraphics[width=6.3cm]{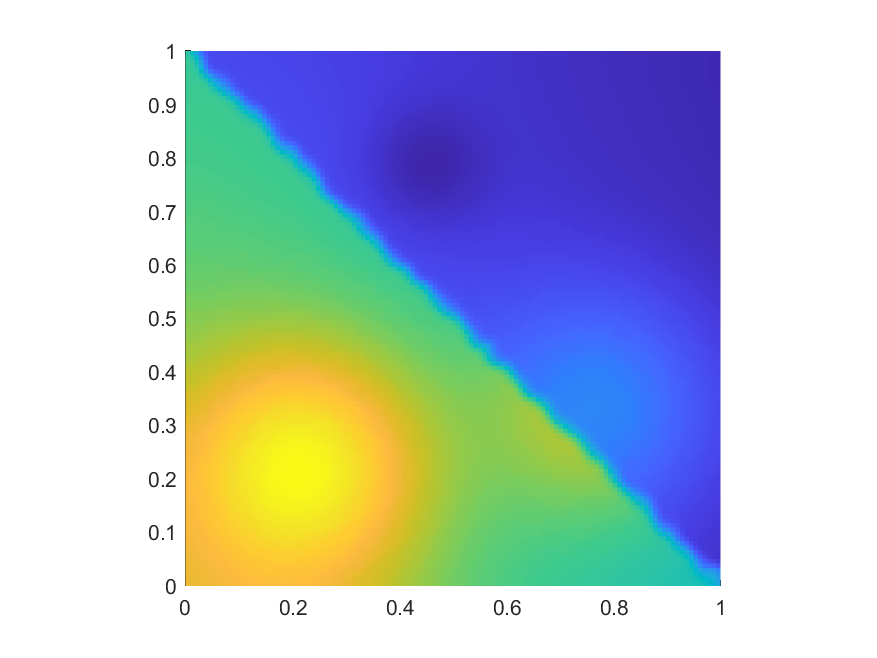} & \hspace{-1cm}	
			\includegraphics[width=6.3cm]{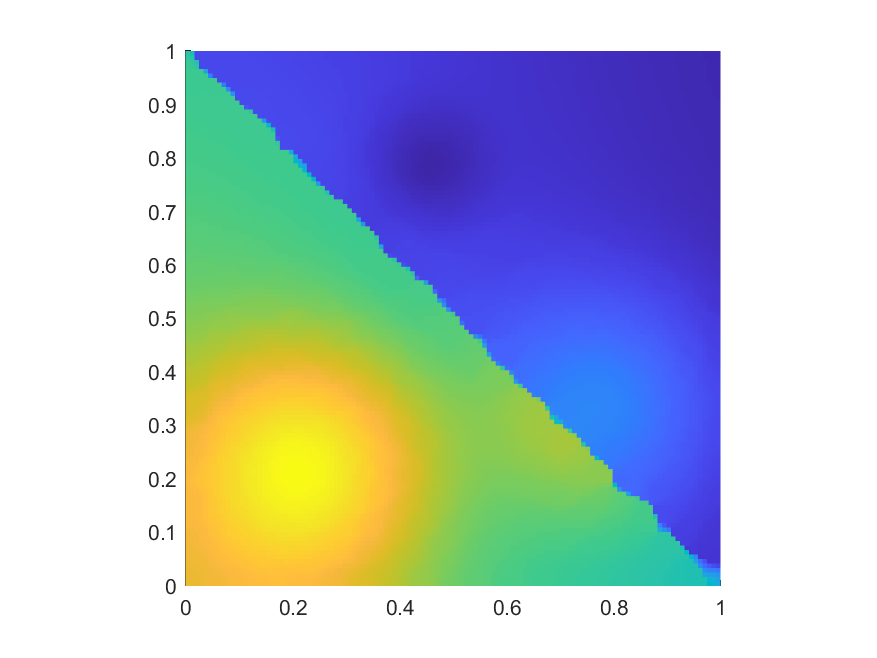} & \hspace{-1cm}
			\includegraphics[width=6.3cm]{diagonaloriginal.eps}\\
\hspace{-1cm}		\includegraphics[width=6.3cm]{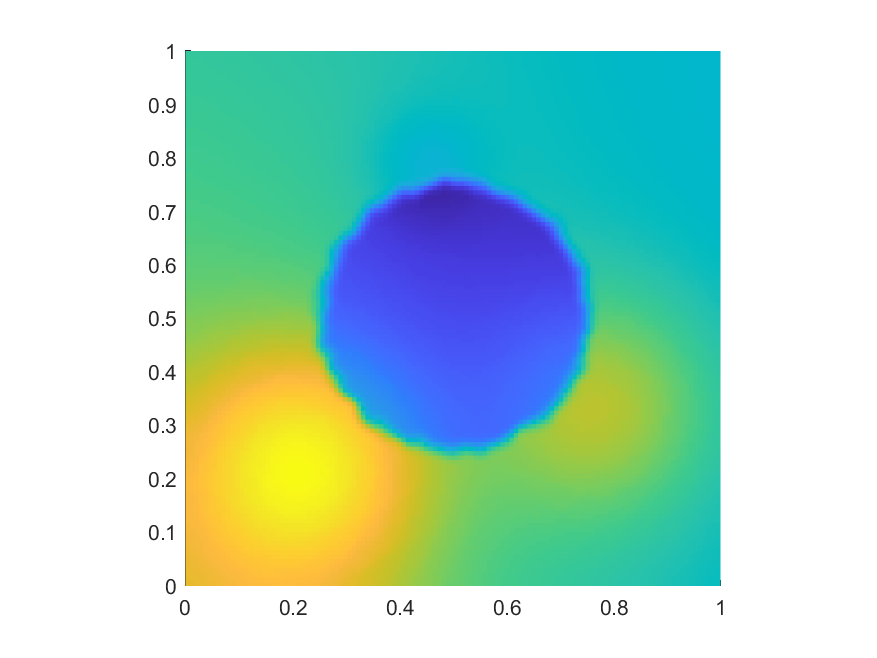} & \hspace{-1cm}	
			\includegraphics[width=6.3cm]{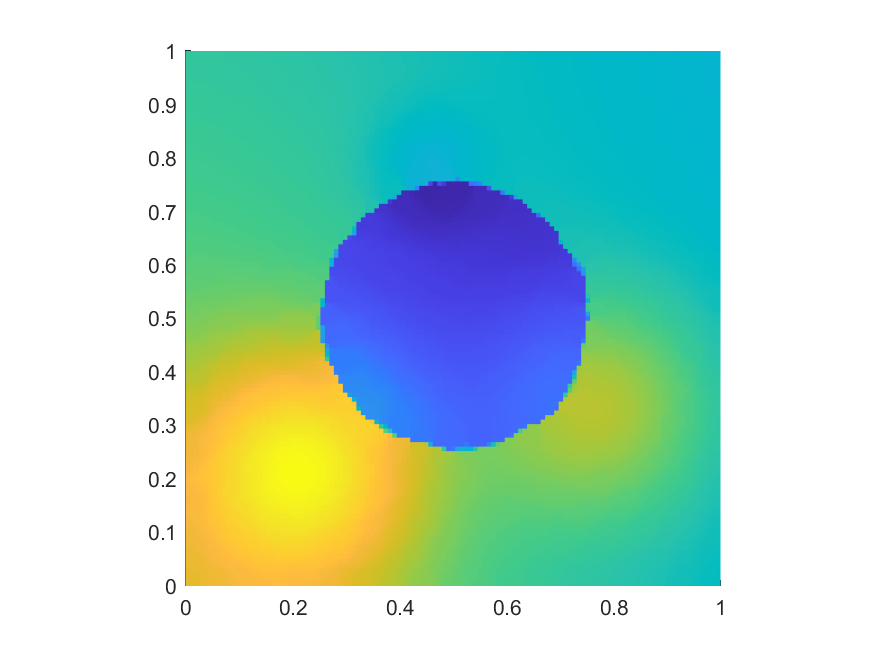} & \hspace{-1cm}
			\includegraphics[width=6.3cm]{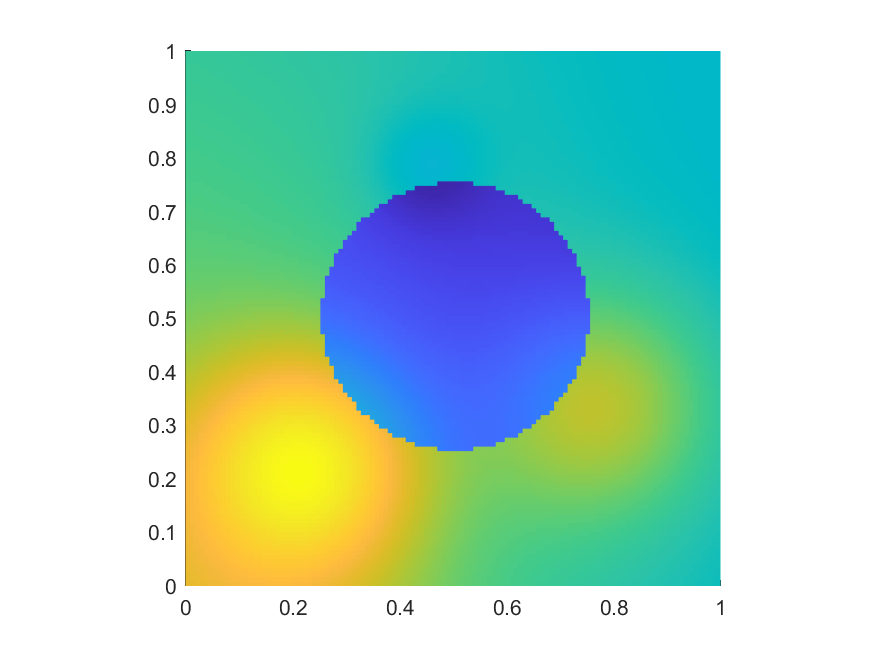}\\
\hspace{-1cm}		\includegraphics[width=6.3cm]{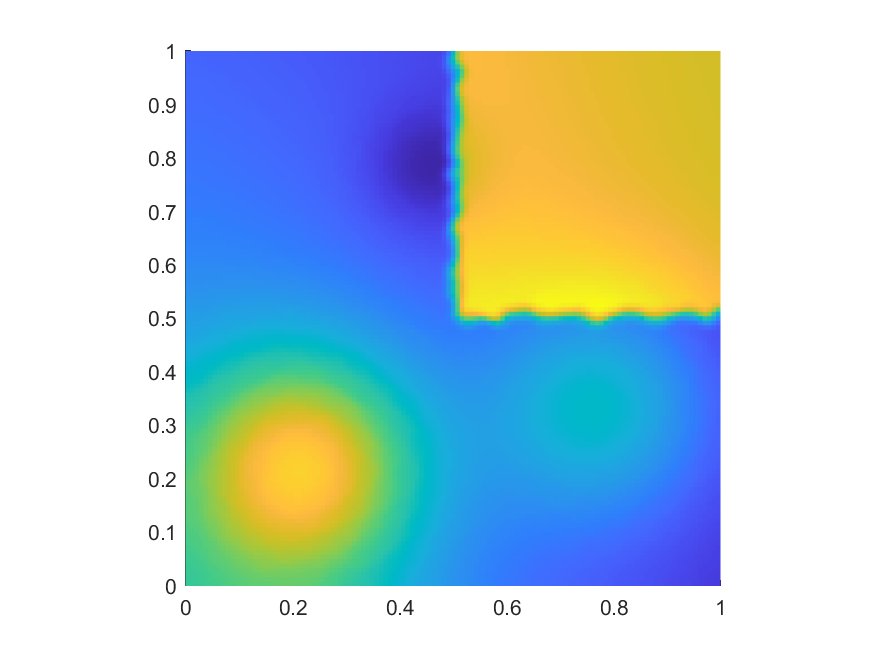} & \hspace{-1cm}	
			\includegraphics[width=6.3cm]{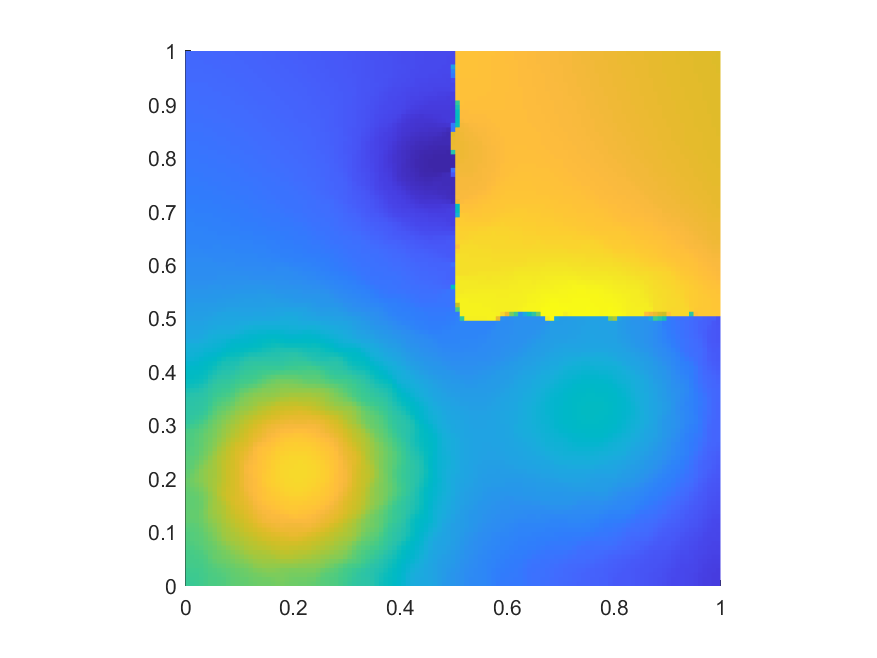} & \hspace{-1cm}
			\includegraphics[width=6.3cm]{cuadradooriginal.eps}
		\end{tabular}
\end{center}
\caption{Approximation to function $\tilde{f}$, Eq. \eqref{funciontildef}, in a 2-d plot with different discontinuity curves using Halton scattered data points.}
	\label{Haltonpointsfig}
	\end{figure}

We consider the domain $[0,1]^2$ and a curve $\Gamma=\{(x,y):\gamma(x,y)=0\}$, where $\gamma:\Omega\to \mathbb{R}$ a continuous function and construct
$\Omega^+=\{(x,y)\in[0,1]^2:\gamma(x,y)\geq 0\}$ and $\Omega^-=[0,1]^2\setminus \Omega^+$ and define
\begin{equation}\label{funciontildef}
\tilde f(x,y)=\left\{
                \begin{array}{ll}
                  1+f(x,y), & (x,y)\in \Omega^+ \\
                  f(x,y), & (x,y)\in \Omega^-
                \end{array}
              \right.
\end{equation}
where $f$ is the Franke's function.
We compute three experiments with $\gamma_1(x,y)=1-x-y$, $\gamma_2(x,y)=0.25^2-x^2-y^2$ and finally with a square with vertices $(0.5,0.5),(0.5,1),(1,0.5),(1,1)$ for Halton scattered regular grid data points. It is clear that in both cases, the non-linear method avoids the diffusion effects close to the discontinuities. When a regular grid is used the algorithm adapts to the discontinuity curve, see Figs. \ref{gridfiguracirc} and \ref{gridfigura}. At the Halton points some irregularities appear, Fig. \ref{Haltonpointsfig}.




\vspace{-0.5cm}

\vspace{-0.5cm}


\end{document}